\documentclass[10pt,twoside]{article}
\makeatletter
\@addtoreset{equation}{section}
\makeatother

\def\pref#1{(\ref{#1})}
\def\dsp{\displaystyle}
\def\Frac#1#2{\frac
{
 {\raise.6ex
 \hbox{$\displaystyle#1$}}
}
{
 {\lower.6ex
 \hbox{$\displaystyle#2$}}
 }
}

\def\FG#1#2#3#4{
{}_2F_1\left(
\begin{array}{c}
\begin{array}{cc} \hskip-10pt#1,{\ } #2 \end{array}\\
\begin{array}{c} \hskip-10pt#3 \end{array}
\end{array}
\hskip-8pt;\,#4
\right)}

\input epsf


\def\begeq{\begin{equation}
\begin{array}{ll}
}
\def\endeq{\end{array}\end{equation}}

\def\bigO{{\cal O}}

\def\tfrac#1#2{{{\lower.6ex
\hbox{$\scriptstyle#1$}}\over 
{\raise.7ex
\hbox{$\scriptstyle#2$}}}}

\begin{document}

\title{The ABC of Hyper Recursions}
\author{Amparo Gil\\
        Departamento de Matem\'aticas, Estad\'{\i}stica y 
        Computaci\'on,\\
        Univ. Cantabria, 39005-Santander, Spain.\\
        \\
    Javier Segura\\
        Departamento de Matem\'aticas, Estad\'{\i}stica y 
        Computaci\'on,\\
        Univ. Cantabria, 39005-Santander, Spain.\\   
    \and
    Nico M. Temme\\
    CWI, P.O. Box 94079, 1090 GB Amsterdam, The Netherlands. \\
     { \small e-mail: {\tt
    amparo.gil@unican.es, 
    javier.segura@unican.es, 
    nicot@cwi.nl}}
    }

\maketitle

\begin{center}
    Dedicated to Roderick Wong on occasion of his $60^{\rm th}$ birthday.
\end{center}

\begin{abstract} 
Each family of Gauss hypergeometric functions 
$$
f_n={}_2F_1(a+\varepsilon_1n, b+\varepsilon_2n ;c+\varepsilon_3n; z),
$$
for fixed $\varepsilon_j=0,\pm1$ (not all $\varepsilon_j$ equal to zero)
satisfies a second order linear difference equation of the form
$$
A_nf_{n-1}+B_nf_n+C_nf_{n+1}=0.
$$
Because of symmetry relations and
functional relations for the Gauss functions, many of the 26 cases (for 
different $\varepsilon_j$ values) can be
transformed into each other.
We give a set of basic equations from which all other equations
can be obtained. For each basic equation, we study the existence
of minimal solutions   
and the character of $f_n$ (minimal or dominant) as $n\to \pm\infty$. 
A second independent solution is given in each basic
case which is dominant when $f_n$ is minimal and vice-versa.
In this way, satisfactory pairs of linearly independent solutions
for each of the 26 second order linear difference equations can be obtained.

\end{abstract}

\par\noindent
Keywords \& Phrases:
Gauss hypergeometric functions,
recursion relations,
difference equations,
stability of recursion relations,
numerical evaluation of special functions,
asymptotic analysis.

\section{Introduction}\label{sec:intro}
The Gauss hypergeometric functions 
\begin{equation}\label{eq:i2}
f_n=\FG{a+\varepsilon_1n}{b+\varepsilon_2n }{c+\varepsilon_3n}{z},
\end{equation}

\noindent where
\begin{equation}\label{eq:i1}
\FG {a}bcz = \sum_{n=0}^\infty \frac{(a)_n (b)_n}{(c)_n\,n!} z^n,\quad |z|<1
\end{equation}
satisfy three-term recurrence relations of the form
\begin{equation}\label{eq:rec}
A_nf_{n-1}+B_nf_n+C_nf_{n+1}=0.
\end{equation}
Examples are given in \cite[p. 558]{Abramowitz:1964:HMF}. 

In this paper we consider recursions with respect to $n$ for the cases
\begin{equation}\label{eq:i3}
\varepsilon_j=0,\pm1, \quad j=1,2,3,\quad 
\varepsilon_1^2+\varepsilon_2^2+\varepsilon_3^2 \ne 0.
\end{equation}
\noindent
and we study the condition of these recurrences by obtaining the regions
in the $z$-complex plane where a minimal solution exists. 

A solution $f_n$ of the recurrence relation \pref{eq:rec}  is said to be
 {\em minimal} if there exists a linearly independent solution $g_n$, of
 the same recurrence relation such that $f_n/g_n\to 0$ as $n\to\infty$. 
In that case $g_n$ is called a {\em dominant} solution. When a recurrence
admits a minimal solution (unique up to a constant factor), this solution
should be included in any numerically satisfactory pair of solutions of
the recurrence. Given a solution of the recurrence, it is crucial to know
the character of the solution (minimal, dominant or none of them) in order
to apply the recurrence relation in a numerically stable way. Indeed, if
$f_n$ is minimal as $n\rightarrow +\infty$, forward recurrence (increasing $n$)
is an ill conditioned process because small initial errors will generally dominate
the recursive solution by introducing an initially small component of a dominant solution; 
backward recurrence is well conditioned in this case. The opposite situation takes place 
for dominant solutions.

The problem of determining the $z$-values for which the $f_n$ functions are minimal or
dominant is considered in this paper. For the different
recursion formulas
we give a companion solution to $f_n$ which, together with $f_n$, form a numerically
satisfactory pair of the corresponding three-term recurrence relation.

\section{Basic recursion relations and their solutions}\label{sec:basic}

There are 26 recursion relations for these functions for all 
choices of $\varepsilon_j$. However, we can use several functional relations in order
to reduce our study to few basic recursions equations. First, we have the symmetry relation
\begin{equation}\label{eq:i4}
\FG abcz=\FG bacz .
\end{equation}
In addition, the following relations can be used
\begin{equation}\label{eq:i5}
\FG {a}bcz = (1-z)^{-a}\,\FG{a}{c-b}{c}{\frac{z}{z-1}}, 
\end{equation}
\begin{equation}\label{eq:i6}
\FG {a}bcz = (1-z)^{-b}\,\FG{c-a}{b}{c}{\frac{z}{z-1}},
\end{equation}
\begin{equation}\label{eq:i7}
\FG {a}bcz = (1-z)^{c-a-b}\,\FG{c-a}{c-b}{c}{z}.
\end{equation}
See \cite[p. 559]{Abramowitz:1964:HMF}. These relations
show that a  great number of the 26 cases are equivalent 
or follow from each other.
In fact we can find 5 basic forms that need to be studied, 
while the remaining cases follow from
the relations in \pref{eq:i4}--\pref{eq:i7}. Few cases for which
$\varepsilon_j=2$ will be treated: only those which follow directly from
the 26 cases for which $|\varepsilon_j|\le 1$, $j=1,2,3$.

When $c=0,-1,-2,\ldots$ the Gauss hypergeometric function in \pref{eq:i1} is not 
defined. Also, when $a$ or $b$ assume non-positive integer values, the series
in \pref{eq:i1} terminates. In the following we will not distinguish about these 
special cases and the general cases, and we assume that all 
representations of the functions to be given are well defined, 
and we will not specify that special values 
of the parameters should be excluded in the results.

In Table \ref{tab:basic} we give an overview of all possible recursions, 
and indicate the five basic forms. Observe that we take recursion in positive direction
equivalent with recursion in negative direction, however we will need 
to distinguish between both directions when studying the asymptotic behavior of the solutions; 
see Section \ref{sec:asy}. Apart
from the notation in Table 1, we will also use the notation 
$(\mbox{sign}(\varepsilon_1)\, \mbox{sign}(\varepsilon_2)\, \mbox{sign}(\varepsilon_2))$ when 
$|\varepsilon_j|\le 1$; for
instance, the case $k=2$ will be also represented as $(++0)$.

\begin{table}
\caption{Only 5 basic forms remain
\label{tab:basic}}
\begin{centering}
{\tt
  \begin{tabular}{||c|c|c|c|l|l||}
 \hline 
   $k$ & $\varepsilon_1$ & $\varepsilon_2$ & $\varepsilon_3$ & type & comments\\ \hline
   1   &   1  &  1   &   1  &  $\equiv 13$ & use \pref{eq:i7}\\ 
   2   &   1  &  1   &   0  &  basic form & \\ 
   3   &   1  &  1   &   -1  &  basic form & \\ 
   4   &   1  &  0   &   1  &  $\equiv 13$ & use \pref{eq:i6}\\ 
   5   &   1  &  0   &   0  &  basic form & \\ 
   6   &   1  &  0   &   -1  &  basic form & \\ 
   7   &   1  &  -1   &   1  &  $\equiv 16$ & use \pref{eq:i6}\\ 
   8   &   1  &  -1   &   0  &  $\equiv 2$ &  use \pref{eq:i5}\\
   9   &   1  &  -1   &   -1  &  $\equiv 6$ & use \pref{eq:i5}\\ \hline
  10   &   0  &  1   &   1  &  $\equiv 13$ & use \pref{eq:i5}\\ 
  11   &   0  &  1   &   0  &  $\equiv 5$ & use \pref{eq:i4}\\ 
  12   &   0  &  1   &   -1  &  $\equiv 6$ & use \pref{eq:i4}\\ 
  13   &   0  &  0   &   1  &  basic form & \\ 
  14   &   0  &  0   &   0  &             & void\\ 
  15   &   0  &  0   &   -1  &  $\equiv 13$ & change sign in 13\\ 
  16   &   0  &  -1   &   1  &  $\equiv 12$ & change signs in 12\\ 
  17   &   0  &  -1   &   0  &  $\equiv 11$ & use \pref{eq:i5}\\ 
  18   &   0  &  -1   &   -1  &  $\equiv 15$ & use \pref{eq:i5}\\ \hline
  19   &   -1  &  1   &    1  &  $\equiv 7$ & use \pref{eq:i4}\\ 
  20   &   -1  &  1   &    0  &  $\equiv 8$ & use \pref{eq:i4}\\ 
  21   &   -1  &  1   &   -1  &  $\equiv 9$ & use \pref{eq:i4}\\ 
  22   &   -1  &  0   &    1  &  $\equiv 16$ & use \pref{eq:i4}\\ 
  23   &   -1  &  0   &    0  &  $\equiv 17$ & use \pref{eq:i4}\\ 
  24   &   -1  &  0   &   -1  &  $\equiv 18$ & use \pref{eq:i4}\\ 
  25   &   -1  &  -1   &   1  &  $\equiv 3$ & change signs in 3\\ 
  26   &   -1  &  -1   &   0  &  $\equiv 2$ & use \pref{eq:i7}\\ 
  27   &   -1  &  -1   &  -1  &  $\equiv 15$ & use \pref{eq:i7}\\ \hline
    \end{tabular}}\\
\end{centering}
\end{table}

\subsection{Selection of a second solution}

Once we have reduced the number of basic recursions to be studied to 5, we will  
give for each basic form the 
coefficients $A_n, B_n$ and $C_n$ of the
recursion relation  \pref{eq:rec} and we will study the character of $f_n$ as a solution
of the corresponding recurrence relation. For this, we will need to find a second solution of the 
relation which forms a satisfactory pair of solutions together with $f_n$ (i.e. a pair which
 includes the minimal solution when it exists).

This second solution is chosen from several connection formulas between the hypergeometric functions.
For example, the functions
\begin{equation}\label{eq:b1}
\FG abcz \quad {\rm and}\quad z^{1-c}\FG{a-c+1}{b-c+1}{2-c}{z}
\end{equation}
satisfy the same differential equation; see \cite[p. 112]{Temme:1996:SFI}.
Both functions in \pref{eq:b1} also satisfy the same 
difference equation, when a suitable normalization 
for the second function is chosen. To obtain this normalization we use the relation
\begin{equation}\label{eq:b2}
\begin{array}{l}
\dsp{\FG abcz}
\dsp{\quad=\frac{\Gamma(c)\Gamma(c-a-b)}{\Gamma(c-a)\Gamma(c-b)}\,\FG{a}{b}{a+b-c+1}{1-z}}\\
\\  \\
\quad\quad\dsp{+\frac{\Gamma(c)\Gamma(a+b-c)}{\Gamma(a)\Gamma(b)}\,(1-z)^{c-a-b}\,\FG{c-a}{c-b}{c-a-b+1}{1-z}}
\end{array}
\end{equation}
(see see \cite[Eq. 15.3.6]{Abramowitz:1964:HMF} or \cite[p. 113]{Temme:1996:SFI}). By replacing $z$ by $1-z$ and $a+b-c+1$ by $c$  
we can write this in the form
\begin{equation}\label{eq:b3}
\begin{array}{l}
\dsp{\FG abcz=P\,\FG{a}{b}{a+b-c+1}{1-z}}\\
\\  \\
\quad\quad\quad\quad\dsp{-\ Q\,z^{1-c}\,\FG{a-c+1}{b-c+1}{2-c}{z},}
\end{array}
\end{equation}
where
\begin{equation}\label{eq:b4}
P=\frac{\Gamma(a+1-c)\,\Gamma(b+1-c)}{\Gamma(a+b-c+1)\,\Gamma(1-c)},\quad 
Q=\frac{\Gamma(c-1)\,\Gamma(a+1-c)\,\Gamma(b+1-c)}{\Gamma(a)\,\Gamma(b)\,\Gamma(1-c)}.
\end{equation}
Another connection formula to be used is
\begin{equation}
\begin{array}{ll}\dsp
&\FG abcz=\\ &\\
&\ \ \Frac{\Gamma(1-a)\Gamma(b-c+1)}{\Gamma(1-c)\Gamma(b-a+1)}
 (-z)^{a-c}(1-z)^{c-a-b}\FG{1-a}{c-a}{b-a+1}{\Frac1{z}} \\&\\
&-\Frac{\Gamma(c-1)\Gamma(b-c+1)\Gamma(1-a)}
{\Gamma(b)\Gamma(c-a)\Gamma(1-c)}
 (-z)^{1-c}(1-z)^{c-a-b}\FG{1-a}{1-b}{2-c}z,
\label{eq:b5}
\end{array}
\end{equation}
where $|\mbox{phase}(-z)|<\pi$. This relation follows from \cite[Eq. 15.3.7]{Abramowitz:1964:HMF},
after we change $z \to 1/z$ and use the result \pref{eq:i7}.

We can take one of the two terms in the  right-hand sides of 
\pref{eq:b2}, \pref{eq:b3} and
\pref{eq:b5} as a second solution of 
the equation that is satisfied by the function in the left-hand side. 
We need to verify if this second solution is linearly independent of the first solution.

When selecting one of these candidates for the second solution, we can skip the gamma 
functions and other terms that are
constant in the recursion.
Also, gamma functions of the form $\Gamma(a-n)$ will be replaced by using the relation
(see \cite[p. 74]{Temme:1996:SFI})
\begin{equation}\label{eq:b6}
\Gamma(a-n)=\frac{(-1)^n \pi}
{\sin\pi a \Gamma(n+1-a)}.
\end{equation}

\subsection{Basic forms and their solutions}

Now we summarize the main results that will be obtained for each basic recursion. We provide
the coefficients of the recursion as well as a satisfactory pair of independent solutions, giving
information on the character of the $f_n$ functions. In later sections we will obtain the regions
of existence of minimal solutions by means of Perron theorem and we will identify the minimal and
a dominant solution by analyzing the asymptotic behaviour of well chosen solutions.

\subsubsection{Basic form {\protect\boldmath $k=2$}}\label{sec:k2}
The $(++0)$ recursion relation reads
\begin{equation}\label{eq:k21}
A_2(a+n, b+n)y_{n-1}+B_2(a+n, b+n)y_n+C_2(a+n, b+n)y_{n+1}=0,
\end{equation}
where
\renewcommand{\arraystretch}{1.5}
\begin{equation}\label{eq:k22}
\begin{array}{ll}
&A_2(a,b)=(c-a)(c-b)(c-a-b-1),\\
&B_2(a,b)=(c-a-b)\{c(a+b-c)+c-2ab +\\
&\quad\quad+z[(a+b)(c-a-b)+2ab+1-c]\},\\
&C_2(a,b)=ab(c-a-b+1)(1-z)2,
\end{array}
\end{equation}
\renewcommand{\arraystretch}{1.0}%
with solutions given by
\begin{equation}\label{eq:k23}
\begin{array}{l}
f_n=\FG{a+n}{b+n}{c}{z},\\ \\ 
g_n=\dsp{\frac{\Gamma(a+n+1-c)\,\Gamma(b+n+1-c)}{\Gamma(a+b+2n-c+1)}}\FG{a+n}{b+n}{a+b+2n-c+1}{1-z}.
\end{array}
\end{equation}
The second solution is taken from the first term in Eq. (\ref{eq:b3}). 

As we will later show, when $z\le 0$ the recurrence has no minimal solutions,
whereas in compact domains
that do not contain points of 
$(-\infty,0]$, $f_n$ is a dominant solution and $g_n$ is minimal.

This case has applications for Jacobi polynomials. We have
\begin{equation}\label{eq:d23}
P_n^{(\alpha,\beta)}(x)={n+\alpha\choose n}
\left(\frac{1+x}{2}\right)^n\FG{-n}{-\beta-n}{\alpha+1}{z},
\quad z=\frac{x-1}{x+1}.
\end{equation} 
A representation with $+n$ at the $a$ and $b$ places follows from
applying \pref{eq:i7}. In the interval of orthogonality $-1\le x\le 1$, we have $z\le0$, and if $x\in[-1,1]$
the recursion relation of the Jacobi polynomials can be used for computing these
functions in forward direction. Only the usual rounding errors should be taken into account.

Notice also how, in addition, information iis obtained regarding the recursion satisfied by 
\begin{equation}
f_n=\FG{a+n}{b+n}{c+2n}{z}. 
\end{equation}
\noindent
By replacing $z$ by $1-z$ and $c$ by $a+b-c$ (see also Eq. (\ref{eq:b6})) we
see that 
\begin{equation}
g_n=\dsp{\frac{\Gamma (c+2n)}{\Gamma (n+1-b+c)\Gamma (n+1-a+c)}}\FG{a+n}{b+n}{a+b-c}{1-z}
\end{equation}
\noindent is also a solution of the same recurrence, $f_n$ being minimal when $z$ is in compact domains 
not containing points in the interval 
$[1,+\infty)$.

\subsubsection{Basic form {\protect\boldmath $k=3$}}\label{sec:k3}
The $(++-)$ recursion relation reads
\begin{equation}\label{eq:k31}
\begin{array}{l}
A_3(a+n, b+n, c-n)y_{n-1}+B_3(a+n, b+n, c-n)y_n\ +\\
\quad\quad\quad C_3(a+n, b+n, c-n)y_{n+1}=0,
\end{array}
\end{equation}
where
\renewcommand{\arraystretch}{1.5}
\begin{equation}\label{eq:k32}
\begin{array}{ll}
&A_3(a,b,c)=-(a-c)(a-c-1)(b-1-c)(b-c)zU,\\
&B_3(a,b,c)=c[c_1U+c_2V+c^3UV],\\
&c_1=(1-z)(b-c)(b-1)[a-1+z(b-c-1)],\\
&c_2=b(b+1-c)(1-z)(a+bz-cz+2z),\\
&c_3=c-2b-(a-b)z,\\
&C_3(a,b,c)=abc(c-1)(1-z)^3V,\\
&U=z(a+b-c+1)(a+b-c+2)+ab(1-z),\cr
&V=(1-z)(1-a-b+ab)+z(a+b-c-1)(a+b-c-2),
\end{array}
\end{equation}
\renewcommand{\arraystretch}{1.0}%
with solutions given by
\renewcommand{\arraystretch}{1.5}
\begin{equation}\label{eq:k33}
\begin{array}{l}
f_n=\FG{a+n}{b+n}{c-n}{z},\\ \\ 
g_n=\dsp{\frac{(-z)^n\Gamma(a+1-c+2n)\,\Gamma(b+1-c+2n)}
{\Gamma(a+n)\Gamma(b+n)\Gamma(2-c+n)\Gamma(1-c+n)}}\ \times \\
\quad\quad\quad\FG{a-c+1+2n}{b-c+1+2n}{2-c+n}{z}.
\label{k=3}
\end{array}
\end{equation}
The second solution is obtained from the second term in Eq. (\ref{eq:b3}) and applying
Eq. (\ref{eq:b6}).

Later we will show that $f_n$ is minimal 
as $n\rightarrow +\infty$ on compact domains inside the curve $r=-9+6\sqrt{3}\cos\frac12\theta$ with
$-\frac13\pi\le \theta\le \frac13\pi$ (Figure 1) and that it is dominant in compact domains outside the curve. 
The opposite situation takes place for $g_n$.

As $n\rightarrow -\infty$, which corresponds to case $k=25$ (i.e., $(-\,-\,+)$), the
roles of $f_n$ and 
$g_n$ are reversed, that is: $f_n$ is dominant (minimal) inside
(outside) the curve represented in Figure 1.

Similarly as happened for the case $k=2$, we can also obtain the condition of the recursion for
a recurrence outside our initial target of 26 cases. From Eq. (\ref{k=3}) we see that the present case is related to the
recursion for $(\varepsilon_1\,\varepsilon_2\,\varepsilon_2)=(2\,2\,1)$. In addition, using Eqs. (\ref{eq:i4}) and
(\ref{eq:i5}) we additionally see that case $(2\,2\,1)$ is related to $(2\,-1\,1)$ and $(-1\,2\,1)$.

\renewcommand{\arraystretch}{1.0}

\subsubsection{Basic form {\protect\boldmath $k=5$}}\label{sec:k5}
The $(+\,0\,0)$ recursion relation reads
\begin{equation}\label{eq:k51}
A_5(a+n)y_{n-1}+B_5(a+n)y_n+C_5(a+n)y_{n+1}=0,
\end{equation}
where
\renewcommand{\arraystretch}{1.5}
\begin{equation}\label{eq:k52}
\begin{array}{ll}
&A_5(a)=(c-a),\\
&B_5(a)=2a-c-(a-b)z\\
&C_5(a)=a(z-1),
\end{array}
\end{equation}
\renewcommand{\arraystretch}{1.0}%
with solutions given by
\renewcommand{\arraystretch}{1.5}
\begin{equation}\label{eq:k53}
\begin{array}{l}
f_n=\FG{a+n}{b}{c}{z},\\
g_n=\dsp{(1-z)^{-n}\frac{\Gamma(a+n+1-c)}{\Gamma(a+n)}\FG{1-a-n}{1-b}{2-c}{z},}
\end{array}
\end{equation}
\renewcommand{\arraystretch}{1.0}%
where we have used the first term in the right-hand side of \pref{eq:b5}.

As we will later show, in compact domains inside the circle $|z-1|=1$, $f_n$ is
dominant and $g_n$ is minimal, while the contrary takes place in compact domains
outside the circle.

\subsubsection{Basic form {\protect\boldmath $k=6$}}\label{sec:k6}
The $(+\,0\,-)$ recursion relation reads
\begin{equation}\label{eq:k61}
A_6(a+n,c-n)y_{n-1}+B_6(a+n,c-n)y_n+C_6(a+n,c-n)y_{n+1}=0,
\end{equation}
where
\renewcommand{\arraystretch}{1.5}
\begin{equation}\label{eq:k62}
\begin{array}{ll}
&A_6(a,c)=z(a-c)(a-c-1)(b-c)[a+z(b+1-c)],\\
&B_6(a,c)=c[a(a-1)(c-1)+a(a-1)(a+3b-4c+2)z\ +\\
&\quad\quad\quad(b-c)(b+1-c)(4a-c-1)z^2-(a-b)(b-c)(b+1-c)z^3],\\
&C_6(a,c)=-ac(c-1)[a-1+z(b-c)](1-z)^2,
\end{array}
\end{equation}
\renewcommand{\arraystretch}{1.0}%
with solutions given by
\renewcommand{\arraystretch}{1.5}
\begin{equation}\label{eq:k63}
\begin{array}{l}
f_n=\FG{a+n}{b}{c-n}{z},\\ \\ 
g_n=\dsp{\frac{\Gamma(a+1-c+2n)\Gamma(b+1-c+n)}
{\Gamma(a+b-c+1+2n)\Gamma(1-c+n)}}\ 
\FG{a+n}{b}{a+b-c+1+2n}{1-z}.
\end{array}
\label{k=6}
\end{equation}
\noindent where $g_n$ is selected from the first term of Eq. (\ref{eq:b3}).

In this case $f_n$ is never minimal. A satisfactory companion of $g_n$ is always
$h_n=g_n-f_n$. The $g_n$ is minimal in compacts domains inside the region between the inner and outer curves in
Figure 2. Contrary, $h_n$ is minimal in compact domains inside the complementary region. 
As $n\rightarrow -\infty$ (case $k=22$, that is $(-\,0\,+)$ or equivalently
 $(0\,-\,+)$, that is, $k=16$), the
role which was played by $h_n$ for the $(+\,0\,-)$ case is now played by $g_n$ and the role of $g_n$ is played by $f_n$.  

From Eq. (\ref{k=6}), we observe that information can also be obtained for the recurrence $(1\,0\,2)$ and then (by 
(\ref{eq:i4}), (\ref{eq:i5}), and (\ref{eq:i7})) 
also for the recurrences $(2\, 1\, 2)$, $(1\, 2\, 2)$, and $(1\, 1\, 2)$.

\renewcommand{\arraystretch}{1.0}%

\subsubsection{Basic form {\protect\boldmath $k=13$}}\label{sec:k13}
The $(0\,0\,+)$ recursion relation reads
\begin{equation}\label{eq:k131}
A_{13}(c+n)y_{n-1}+B_{13}(c+n)y_n+C_{13}(c+n)y_{n+1}=0,
\end{equation}
where
\renewcommand{\arraystretch}{1.5}
\begin{equation}\label{eq:k132}
\begin{array}{ll}
&A_{13}(c)=c(c-1)(z-1),\\
&B_{13}(c)=c[c-1-(2c-a-b-1)z],\\
&C_{13}(c)=(c-a)(c-b)z,
\end{array}
\end{equation}
\renewcommand{\arraystretch}{1.0}%
with solutions given by
\begin{equation}\label{eq:k133}
\begin{array}{l}
f_n=\FG{a}{b}{c+n}{z},\\ \\ 
g_n=\dsp{\frac{(-1)^n(1-z)^n\Gamma(c+n)}
{\Gamma(c-a-b+1+n)}}\ \FG{c-a+n}{c-b+n}{c-a-b+1+n}{1-z},
\end{array}
\end{equation}
where we have used the second term in the right-hand side of \pref{eq:b2}.

As we will show next, $f_n$ is minimal in compact domains inside the region $\Re z<1/2$, where
$g_n$ is dominant. Contrary, in compact domains inside the region $\Re z>1/2$ $f_n$ is dominant and
$g_n$ is minimal.

When $n\rightarrow -\infty$ (that is, for the case $k=15$, corresponding to $(0\,0\,-)$) $f_n$ is
a dominant solution both for $\Re z<1/2$ and $\Re z>1/2$. The satisfactory companion solution 
for $f_n$ (the minimal solution) is chosen differently for the cases $\Re z <1/2$ and $\Re z>1/2$
(see \S \ref{sec:k15}).

\section{Domains for minimal and dominant solutions}\label{sec:dom}

Perron's theorem (see \cite[Appendix B]{Wimp:1984:CRR}) gives in the case of finite limits
the following results.
Let $\alpha$ and $\beta$ denote
\begin{equation}\label{eq:d1}
\alpha:=\lim_{n\to\infty}\frac{B_n}{C_n},\quad \beta:=\lim_{n\to\infty}\frac{A_n}{C_n}.
\end{equation} 

Let $t_1$ and $t_2$ denote the zeros of the characteristic polynomial $t^2+\alpha t+\beta=0$.
If $|t_1|\ne |t_2|$, then the difference equation \pref{eq:rec} has two 
linear independent solutions $f_n$ and $g_n$
with the properties
\begin{equation}\label{eq:d2}
\frac{f_{n+1}}{f_n}\sim t_1,\quad \frac{g_{n+1}}{g_n}\sim t_2.
\end{equation} 
If $|t_1| = |t_2|$, then 
\begin{equation}\label{eq:d3}
\limsup_{n\to\infty}|y_n|^{\frac1n}=|t_1|
\end{equation} 
for any non-trivial solution $y_n$ of \pref{eq:rec}.

In the following subsections, we give for the five basic forms the domains in the $z-$plane where 
$|t_1|\ne|t_2|$. In these domains there is a true distinction between the two solutions of
\pref{eq:rec}. If $|t_1|>|t_2|$ then the solution $f_n$ that satisfies 
the relation in \pref{eq:d2}
is a maximal solution and $g_n$ is the minimal solution.
On the curves where $|t_1|=|t_2|$ the two solutions are neither dominant nor minimal, and 
recursion in forward or backward direction is not unstable.

For all five basic forms the ratios $A_n/C_n$ and $B_n/C_n$ of the difference equation
\pref{eq:rec} tend to finite limits as $|n|\to \infty$. Interestingly, in all cases these limits 
are functions of $z$, and they
are not depending on the parameters $a$, $b$ or $c$.

\subsection{The domains for basic form {\protect\boldmath $k=2$}}\label{sec:dk2}
The limits $\alpha$ and $\beta$ of \pref{eq:k21} are
\begin{equation}\label{eq:d21}
\alpha=-\frac{2(z+1)}{(1-z)^2}, \quad \beta= \frac{1}{(1-z)^2}.
\end{equation} 
The zeros of the characteristic polynomial are
\begin{equation}\label{eq:d22}
t_1= \frac{1}{(1-\sqrt{z})^2}, \quad t_2= \frac{1}{(1+\sqrt{z})^2}.
\end{equation} 
The equation $|t_1|=|t_2|$ holds when $z\le 0$, otherwise $|t_1|>|t_2|$.

\subsection{The domains for basic form {\protect\boldmath $k=3$}}\label{sec:dk3}
\begin{figure}
\begin{center}
\epsfxsize=8cm \epsfbox{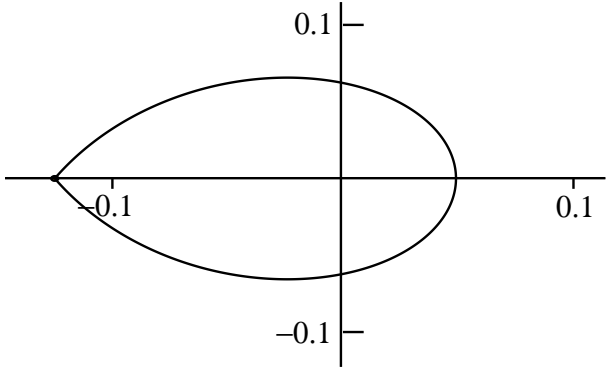}
\end{center}
\caption{\small The curve $|t_1|=|t_2|$ for the basic form $k=3$.\label{fig:abc3}}
\end{figure}

The limits $\alpha$ and $\beta$ of \pref{eq:k31} are
\begin{equation}\label{eq:d31}
\alpha= \frac{8z^2+20z-1}{(1-z)^3}, \quad \beta= -\frac{16z}{(1-z)^3}.
\end{equation} 
The zeros of the characteristic polynomial are
\begin{equation}\label{eq:d32}
t_1= \frac{1-20z-8z^2+(8z+1)^{\frac32}}{2(1-{z})^3}, \quad 
t_2= \frac{1-20z-8z^2-(8z+1)^{\frac32}}{2(1-{z})^3}.
\end{equation} 
We write this in the form
\renewcommand{\arraystretch}{1.5}
\begin{equation}\label{eq:d33}
\begin{array}{ll}
&\dsp{t_1= \frac{27-18w^2-w^4+8w^3}{16(1-{z})^3}=\frac{32(1+{w})}{(3+{w})^3},} \\
&\dsp{t_2= \frac{27-18w^2-w^4-8w^3}{16(1-{z})^3}=\frac{32(1-{w})}{(3-{w})^3},}
\end{array}
\end{equation} 
\renewcommand{\arraystretch}{1.0}
where $w=\sqrt{8z+1}$. To find the curve in the $w-$plane defined by $|t_1|=|t_2|$, 
we write $w=re^{i\theta}$. This gives the curve described by
\begin{equation}\label{eq:d34}
r=-9+6\sqrt{3}\cos\tfrac12\theta,\quad -\tfrac13\pi \le \theta\le \tfrac13\pi.
\end{equation} 
In Figure \ref{fig:abc3} we show this curve in the $z-$plane. 
In the domain interior to this curve we have $|t_1|>|t_2|$.

\subsection{The domains for basic form {\protect\boldmath $k=5$}}\label{sec:dk5}
The limits $\alpha$ and $\beta$ of \pref{eq:k51} are
\begin{equation}\label{eq:d51}
\alpha=\frac{z-2}{1-z}, \quad \beta= \frac{1}{1-z}.
\end{equation} 
The zeros of the characteristic polynomial are
\begin{equation}\label{eq:d52}
t_1= 1, \quad t_2= \frac{1}{1-z}.
\end{equation} 
The equation $|t_1|=|t_2|$ holds when $|1-z|=1$, which defines a circle with centre $z=1$ and radius 1.
Inside the circle we have $|t_2|>|t_1|$.

\begin{figure}
\begin{center}
\epsfxsize=8cm \epsfbox{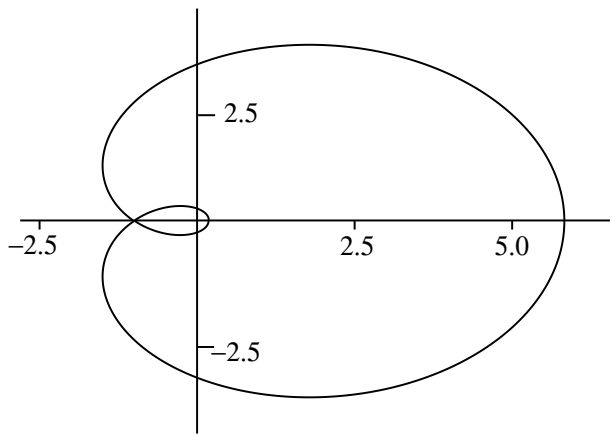}
\end{center}
\caption{\small The curve $|t_1|=|t_2|$ for the basic form $k=6$.\label{fig:abc6}}
\end{figure}

\subsection{The domains for basic form {\protect\boldmath $k=6$}}\label{sec:dk6}
The limits $\alpha$ and $\beta$ of \pref{eq:k61} are
\begin{equation}\label{eq:d61}
\alpha= -\frac{z^2-6z+1}{(1-z)^2}, \quad \beta= -\frac{4z}{(1-z)^2}.
\end{equation} 
The zeros of the characteristic polynomial are
\begin{equation}\label{eq:d62}
t_1= 1, \quad 
t_2= -\frac{4z}{(1-z)^2}.
\end{equation} 
To find the curve  defined by $|t_1|=|t_2|$, 
we write $z=re^{i\theta}$. This gives the curve described by
\begin{equation}\label{eq:d63}
r=2+\cos\theta\pm\sqrt{\cos2\theta+4\cos\theta+3},\quad -\pi \le \theta\le \pi.
\end{equation} 
Both signs give a closed loop with common point $-1$.
In Figure \ref{fig:abc6} we show this curve in the $z-$plane. 
In the domain interior to the inner curve we have $|t_1|>|t_2|$; 
between the inner curve and the outer curve
we have $|t_1|<|t_2|$, and outside the outer curve $|t_1|>|t_2|$.

\subsection{The domains for basic form {\protect\boldmath $k=13$}}\label{sec:dk13}
The limits $\alpha$ and $\beta$ of \pref{eq:k131} are
\begin{equation}\label{eq:d131}
\alpha=-\frac{2z-1}{z}, \quad \beta= \frac{z-1}{z}.
\end{equation} 
The zeros of the characteristic polynomial are
\begin{equation}\label{eq:d132}
t_1= 1, \quad t_2= \frac{z-1}{z}.
\end{equation} 
The equation $|t_1|=|t_2|$ holds when $\Re z =\frac12$.
When $\Re z>\frac12$ we have  $|t_1|>|t_2|$.

\section{Asymptotics for minimal and dominant solutions}\label{sec:asy}
In some cases we simply use the power series in \pref{eq:i1},
which provides an asymptotic expansion for large $c$. In some other cases we can use 
connection formulas for transforming the Gauss function to the case for large $c$.
We also use the integral representations
\begin{equation}\label{eq:as1}
\FG abcz 
=\frac{\Gamma(c)}{\Gamma(b)\Gamma(c-b)}\int_0^1 t^{b-1}(1-t)^{c-b-1}(1-tz)^{-a}\,dt,
\end{equation}
\begin{equation}\label{eq:as2}
\FG abcz =\frac{\Gamma(c)\,\Gamma(1+b-c)}{2\pi i\,\Gamma(b)}\int_0^{(1^+)} t^{b-1}(t-1)^{c-b-1}(1-tz)^{-a}\,dt,
\end{equation} 
where in the first integral $\Re c>\Re b>0$ and in the second one $\Re b>0$.
In the second integral the contour starts and terminates 
at $t=0$ and encircles the point $t=1$ in the positive direction. 
The point $t=1/z$ should be outside the contour. The many-valued functions of 
the integrand assume their principal branches: the phase of $(1-tz)$ 
tends to zero when $z\to0$, and the phases of 
$t$ and $(t-1)$ are zero at the point where the contour cuts the real 
positive axis (at the right of $t=1$).
To prove \pref{eq:as2} integrate (if $\Re c>\Re b>0$) along the interval $(0,1)$ 
with proper choices of the branches of the many-valued functions at the upper and lower 
sides of $(0,1)$, and obtain \pref{eq:as1}; 
see also \cite[p. 111]{Temme:1996:SFI}. 
The integral in \pref{eq:as1} is Euler's
well-known standard representation.

We use saddle-point methods (see \cite{Wong:2001:AAI}) to obtain asymptotic estimates
of these integrals when one of more parameters are large.
We omit details of the saddle-point analysis because we only need to obtain the main terms 
in the asymptotic estimates for identifying minimal and dominant solutions.
We assume that $z$ is fixed and properly inside the domains described in Section~\ref{sec:dom}.

For recent papers on uniform asymptotic expansions of 
hypergeometric functions, see 
\cite{Jones:2001:AHF}, \cite{Olde:2003:UAI} and \cite {Olde:2003:UII}.

\subsection{Asymptotics for basic form {\protect\boldmath $k=2$}}\label{sec:ak2}
For $f_n$ of \pref{eq:k23} we use \pref{eq:as2} and we obtain
\begin{equation}\label{eq:ask21}
f_n =\frac{\Gamma(c)\,\Gamma(1+b+n-c)}{2\pi i\,\Gamma(b+n)}
\int_0^{(1^+)} t^{b-1}(t-1)^{c-b-1}(1-tz)^{-a}\,e^{n\phi(t)}\,dt,
\end{equation} 
where $\phi(t)= \ln t-\ln(t-1)-\ln(1-tz)$. The saddle-points 
are found by putting $\phi'(t)=0$, giving $t_{\pm}=\pm1/\sqrt{z}$. If $z\in(0,1)$
we have $t_+>1$ and
the saddle-point contour starts at $t=0$, turns around $t=1$ through $t_+$, and returns to $t=0$.
The dominant term in the asymptotic estimate is $e^{\phi(t_+)}=t_1^n$ (see \pref{eq:d22}).
This also holds for complex values of $z$.

For $g_n$ of \pref{eq:k23} we use \pref{eq:as1} and we obtain
\renewcommand{\arraystretch}{1.75}
\begin{equation}\label{eq:ask22}
\begin{array}{l}
\dsp{g_n 
=\frac{\Gamma(a+n+1-c)\,\Gamma(b+n+1-c)}
{\Gamma(a+n+1-c)\Gamma(b+n)}}\  \times\\
\quad\quad\quad\dsp{\int_0^1 t^{b-1}(1-t)^{c-b-1}(1-t(1-z))^{-a}\,e^{n\psi(t)}\,dt,}
\end{array}
\end{equation}
\renewcommand{\arraystretch}{1.0}
where $\psi(t)=\ln t+\ln(1-t)-\ln(1-t(1-z))$. 
The saddle-points are $t_+=\sqrt{t_2}$, $t_-=\sqrt{t_1}$. 
When $z\in(0,1)$ we have $t_+\in(0,1)$ and  this point gives the 
dominant contribution. 
The dominant term in the asymptotic estimate is $e^{\psi(t_+)}=t_2^n$.
This also holds for complex values of $z$.

When $z$ is not a negative real number we infer that 
for $g_n$  the main term in the asymptotics is $t_2^n$ and that for
$f_n$ the main term is $t_1^n$. We conclude that in compact domains  that do not contain
points of $(-\infty,0]$, $g_n$ of \pref{eq:k23} is the minimal solution and $f_n$ is a dominant solution of
\pref{eq:k21}.

\subsection{Asymptotics for basic form {\protect\boldmath $k=3$}}\label{sec:ak3}
For $f_n$ of \pref{eq:k33} we use \pref{eq:as2} and we obtain
\begin{equation}\label{eq:ask31}
f_n =\frac{\Gamma(c-n)\,\Gamma(1+b-c+2n)}{2\pi i\,\Gamma(b+n)}
\int_0^{(1^+)} t^{b-1}(t-1)^{c-b-1}(1-tz)^{-a}\,e^{n\phi(t)}\,dt,
\end{equation} 
where $\phi(t)= \ln t-2\ln(t-1)-\ln(1-tz)$. The saddle-points 
are  $t_{\pm}=(1\pm w)/(4z)$, $w=\sqrt{1+8z}$. If $z\in(0,1)$
we have $t_+>1$ and
the saddle-point contour starts at $t=0$, turns around $t=1$ through $t_+$, and returns to $t=0$.
The dominant term in the asymptotic estimate is $e^{\phi(t_+)}=(-1)^n2^{-2n}t_2^n$ 
(see \pref{eq:d32}).
This also holds for complex values of $z$.

For $g_n$ of \pref{eq:k33} we again use \pref{eq:as2} and we obtain
\renewcommand{\arraystretch}{1.75}
\begin{equation}\label{eq:ask32}
\begin{array}{l}
\dsp{g_n 
=\frac{(-z)^n\,\Gamma(a+1-c+2n)}
{ 2\pi i\,\Gamma(a+n)\Gamma(1-c+n)}}\  \times\\
\quad\quad\quad\dsp{\int_0^{(1+)} t^{b-c}(t-1)^{-b}(1-tz)^{-a+c-1}\,e^{n\psi(t)}\,dt,}
\end{array}
\end{equation}
\renewcommand{\arraystretch}{1.0}
where $\psi(t)=2\ln t-\ln(t-1)-2\ln(1-tz)$. 
The saddle-points are $t_{\pm}=(-1\pm w)/(2z)$.
When $z\in(0,1)$ we have $t_+\in(1,2)$ and  the saddle-point contour is similar as for $f_n$.
The dominant term in the asymptotic estimate is $e^{\psi(t_+)}=(-z)^{-n}2^{-2n}t_1^n$
(see \pref{eq:d32}).
This also holds for complex values of~$z$.

Taking into account the asymptotics of the gamma function in 
front of the integrals in \pref{eq:ask31} and \pref{eq:ask32}, 
we infer that the dominant terms in the asymptotic estimate of $f_n$ is $t_2^n$; and 
for $g_n$ it is $t_1^n$.
We conclude that in compact domains  interior to the curve of Figure \ref{fig:abc3}, 
$g_n$ of \pref{eq:k33} is a dominant solution and $f_n$ is the minimal solution of
\pref{eq:k31}. In compact domains exterior to 
this curve the roles of $f_n$ and $g_n$ are interchanged.
 
\subsection{Asymptotics for basic form {\protect\boldmath $k=5$}}\label{sec:ak5}
For $f_n$ of \pref{eq:k53} we use \pref{eq:as1} and we obtain
\begin{equation}\label{eq:ask51}
f_n =\frac{\Gamma(c)}{\Gamma(b)\,\Gamma(c-b)}
\int_0^{1} t^{b-1}(1-t)^{c-b-1}(1-tz)^{-a-n}\,dt,
\end{equation} 
If $|1-z|<1$ the main contributions to the integral come from $t-$values near the end-point
$1$. In that case the dominant term in the asymptotic estimate is $t_2^n$ 
(see \pref{eq:d52}). If $|1-z|>1$ the main contributions to the integral come from 
$t-$values near the end-point
$0$. In that case $f_n=\bigO(1/n)$ as $n\to\infty$.

For $g_n$ of \pref{eq:k53} we use \pref{eq:i6}, which gives
\begin{equation}\label{eq:ask52}
g_n=(1-z)^{1+b-c-n}\frac{\Gamma(a+n+1-c)}{\Gamma(a+n)}\FG{1+a-c+n}{1-b}{2-c}{\frac{z}{1-z}}.
\end{equation} 
We can use the same method as for $f_{n}$, and conclude that 
in compact domains  inside the disk  $|z-1|<1$,
$g_n$ of \pref{eq:k53} is a minimal solution and $f_n$ is a dominant solution of
\pref{eq:k51}, respectively corresponding with $t_1$ and $t_2$ of \pref{eq:d52}. 
In compact domains outside the disk the roles of $f_n$ and $g_{n}$ are
interchanged.

\subsection{Asymptotics for basic form {\protect\boldmath $k=6$}}\label{sec:ak6}
For $g_n$ of \pref{eq:k63} we use \pref{eq:i4} and \pref{eq:as1}, and obtain
\begin{equation}\label{eq:ask61}
g_n=\frac{\Gamma(a+1-c+2n)}{\Gamma(a+n)\Gamma(1-c+n)}\,
\int_0^1 t^{a-1}(1-t)^{b-c}[1-(1-z)t]^{-b} e^{n\phi(t)}\,dt
\end{equation}
where $\phi(t)=\ln t+\ln(1-t)$. The point $t=\frac12$ gives the main contribution,
and the dominant term of the integral is $4^{-n}$. 

Next we consider $h_n:=g_n-f_n$. From \pref{eq:k133}, \pref{eq:b3} and \pref{eq:as1}
it follows that
\renewcommand{\arraystretch}{1.75}
\begin{equation}\label{eq:ask62}
\begin{array}{l}
\dsp{h_n 
=\frac{z^{1-c+n}\,\Gamma(c-n-1)\Gamma(a+1-c+2n)\Gamma(2-c+n)}
{\Gamma(a+n)\Gamma(b)\Gamma(1-c+n)\Gamma(1-b)}}\  \times\\ 
\quad\quad\quad\dsp{\int_0^1 t^{b-c}(1-t)^{-b}[1-t(1-z)]^{-a+c-1}\,e^{n\psi(t)}\,dt,}
\end{array}
\end{equation}
\renewcommand{\arraystretch}{1.0}%
where $\psi(t)=\ln t-2\ln(1-tz)$. The end-point $t=1$ gives the dominant contribution
$z^n/(1-z)^{2n}$ to the integral. 

Taking into account the contribution from the gamma functions
in \pref{eq:ask61} and \pref{eq:ask62}, it follows that $g_n$ corresponds 
with $t_1$ of \pref{eq:d62}, for all $z$, and $h_n$ with $t_2$.

Because $f_n=g_n-h_n$, this function is never a minimal solution of
\pref{eq:k61}. 
In compact domains interior to the inner  curve and outside the outer curve
(see  Figure \ref{fig:abc6}), $g_n$ is a dominant solution  and 
$h_n$ the minimal solution; in compact domains between the inner and the outer 
curve the roles of $g_n$ and $h_n$ are interchanged.

\subsection{Asymptotics for basic form {\protect\boldmath $k=13$}}\label{sec:ak13}
For $f_n$  of \pref{eq:k133} we have the estimate $f_n=1+\bigO(1/n)$ for all $z$. 
For $g_n$ we apply \pref{eq:i7}, and obtain
\begin{equation}\label{eq:ask133}
g_n=\frac{(1-z)^n}{z^n}\frac{z^{1-c}\,\Gamma(c+n)}
{\Gamma(c-a-b+1+n)}\ \FG{1-b}{1-a}{c-a-b+1+n}{1-z}.
\end{equation}
We conclude that in compact domains of $\Re z < \frac12$, $f_n$ is the minimal solution 
and $g_n$ is a dominant solution.  In compact domains of  $\Re z > \frac12$ the roles of 
$f_n$ and $g_n$ are interchanged.

\section{The cases {\protect\boldmath $k=15, k=16, k=25$}}\label{sec:chcases}
The cases $k=15, k=16, k=25$ in Table \ref{tab:basic} are special because we refer 
to these as \lq change signs in other cases\rq,
and we do not consider them as basic forms. The recursion relations for the cases
$k=15, k=16, k=25$ are the same as those for $k=13, k=12, k=3$, respectively, 
when we recur backwards, that is to $-\infty$. 
The zeros $t_1$ and $t_2$ of the characteristic polynomial
do not change when we change the recursion direction.
In this section we give linearly independent satisfactory solutions for the cases
$k=15, k=16, k=25$.

\subsection{The case {\protect\boldmath $k=15$}}\label{sec:k15}
We verify if the solutions 
given for the  recursion relation of the case $k=13$ can be used as satisfactory linear 
independent solutions for the same recursion relation used in backward direction.

The asymptotics for of $f_n$ of \pref{eq:k133} for large positive $n$ 
follows easily from the power series. For negative $n$ we can use
\pref{eq:b3} and for the term with $Q$ we use also \pref{eq:i7}. This gives 
\begin{equation}\label{eq:b3mod}
\begin{array}{l}
\dsp{f_{-n}=\FG ab{c-n}z=P\,\FG{a}{b}{a+b-c+n+1}{1-z}}\\
\\  \\
\quad\quad\quad\quad\dsp{-\ Q\,z^{1-c+n}(1-z)^{c-n-a-b}\,\FG{1-a}{1-b}{2-c+n}{z},}
\end{array}
\end{equation}
where $P$ and $Q$ are as in \pref{eq:b4} with $n$ replaced with $-n$.
Again the power series can be used for obtaining the asymptotics, and 
it follows that for $f_{-n}$ of \pref{eq:b3mod} we have
$f_{-n}=1+\bigO(1/n)$ as $n\to+\infty$, when $|z/(1-z)| <1$; 
otherwise, the dominant term in the asymptotic behaviour is $|z/(1-z)|^n$.

For the second solution $g_n$ of  
\pref{eq:k133} we obtain, after using \pref{eq:i7},
\begin{equation}\label{eq:gnk15}
g_{-n}=\frac{(-1)^nz^{1-c}\Gamma(c-n)}{\Gamma(c-a-b+1-n)} 
\left(\frac{z}{1-z}\right)^n
\FG{1-b}{1-a}{c-a-b+1-n}{1-z}.
\end{equation}
We can verify as for the first $F-$function in \pref{eq:b3mod} that the $F-$function 
in \pref{eq:gnk15} has dominant term $[(1-z)/z]^n$ when $|z/(1-z)| <1$; 
otherwise, the dominant term in the asymptotic behaviour is algebraic in $n$. 
Because of the extra factor $[z/(1-z)]^n$ in front of the $F-$function in
\pref{eq:gnk15}, it follows that $f_{-n}$ and $g_{-n}$
are minimal or dominant in the same $z-$domain.
Hence the $\{f_n,g_n\}$ of \pref{eq:k133} does not constitute a 
satisfactory pair of linearly independent solutions for the recursion \pref{eq:k132} in backward direction.
As a consequence, we cannot use case $k=13$ in backward direction 
for the present case $k=15$.

The second solution $g_n$ of  \pref{eq:k133} has been selected by considering the 
second term in \pref{eq:b2}. When instead we consider the first term, 
using \pref{eq:b6}, 
with $c$ replaced with
$c-n$, we can take as another solution of \pref{eq:k132} the function $h_n$,
where
\begin{equation}\label{eq:hnk15}
h_{-n}=\frac {\Gamma(n+1-c+a)\Gamma(n+1-c+b)}{\Gamma(n+1-c)\Gamma(n+1-c+a+b)}
\FG{a}{b}{n+1-c+a+b}{1-z}.
\end{equation}
The $F-$function is $1+\bigO(1/n)$ for large positive $n$ in compact 
$z-$domains. Hence, $h_n$ can be used as a proper second solution
together with $f_{n}$ for the recursion relation of the case
$k=13$ in backward direction when 
$|z/(1-z)| >1$. We still need another function for the complementary domain.
We consider for this the connection formula \pref{eq:b5}
use the first term in the right-hand side, and \pref{eq:b6}, and we take as
a new solution of  \pref{eq:k132} the function $j_n$, where
\begin{equation}
\begin{array}{l}\dsp
{j_{-n}=(-1)^n\Frac{\Gamma(a-c+n+1)\Gamma(b-c+n+1)}
{\Gamma(2-c+n)\Gamma(1-c+n)}} \times \\
\\
\quad\quad\quad \dsp{\left(\frac{z}{1-z}\right)^n\FG{1-a}{1-b}{2-c+n}{z}}.
\label{eq:jnk15}
\end{array}
\end{equation}
The $F-$function is $1+\bigO(1/n)$ for large positive $n$ in compact 
$z-$domains. Hence, $j_n$ can be used as a proper second solution
together with $f_{n}$ for the recursion relation of the case
$k=13$ in backward direction when 
$|z/(1-z)| <1$. 

Summarizing, for the present case $k=15$ we can use the recursion \pref{eq:k132}
of basic form $k=13$ in backward direction with 
$f_n$ of \pref{eq:k133} and $h_n$ of \pref{eq:hnk15} for $z-$values in 
compact domains of the half-plane $\Re z> \frac12$ where $f_n$ is a dominant and $h_n$
is the minimal solution. In compact domains 
of $\Re z<\frac12$ we can use $f_n$ together with $j_n$ of \pref{eq:jnk15}, 
with $f_n$ as a dominant and $j_n$ as the minimal solution.

\subsection{The case {\protect\boldmath $k=16$}}\label{sec:k16}
We verify if the solutions 
given for the  recursion relation of the case $k=6$ can be used as linear 
independent solutions for the same recursion relation used in backward direction.

The function $f_n$ of the case $k=6$ given in \pref{eq:k63} is $1+\bigO(1/n)$ as 
$n\to-\infty$. This easily follows when we use the integral representation in
\pref{eq:as1}. Using this integral
$g_n$ of \pref{eq:k63} it follows that for $g_n$
the dominant term is $t_2^{-n}$, where $t_2$ is given  in 
\pref{eq:d62}.

It follows that for the present case $k=16$ we can use the recursion of basic form $k=6$,
that is,
\pref{eq:k61}, in backward direction with $f_n$ and $g_n$ of \pref{eq:k63}.
Inside compact domains interior to the inner curve (around the origin) 
of Figure \ref{fig:abc6} 
and exterior to the outer curve
$f_n$ is a dominant and $g_n$ is the minimal solution. In compact domains between 
the two curves the roles of
$f_n$ and $g_n$ are interchanged.

\subsection{The case {\protect\boldmath $k=25$}}\label{sec:k25}
We verify if the solutions 
given for the  recursion relation of the case $k=3$ can be used as linear 
independent solutions for the same recursion relation used in backward direction.

When we apply \pref{eq:i7} to both $f_n$ and $g_n$ of \pref{eq:k33},
we see that the asymptotics for $n\to-\infty$ follows from that of $n\to+\infty$,
when the roles of $f_n$ and $g_n$ are interchanged.

We conclude that in compact domains  interior to the curve of Figure \ref{fig:abc3}, 
$f_n$ is a dominant solution and $g_n$ is the minimal solution of
of \pref{eq:k31} when used in backward direction.
In compact domains exterior to 
this curve the roles of $f_n$ and $g_n$ are interchanged.

 \section{Numerical examples}\label{sec:num}
The power series \pref{eq:i1} is very useful for numerical computations for $z$ properly
inside the unit disk. Transformations and connection formulas as in
\pref{eq:i5}, \pref{eq:i6}, \pref{eq:b2} and \pref{eq:b5} can be used 
to cover large parts of the complex
$z-$plane. Other connection formulas are available with other $z-$arguments
and in fact for the computation of the Gauss function we can 
use power series with powers of
\begin{equation}\label{eq:num1}
z,\quad 1-z,\quad \frac1z,\quad \frac{z-1}{z},\quad \frac1{1-z},
\quad \frac{z}{z-1}.
\end{equation}
For numerical computations we need convergence conditions like
\begin{equation}\label{eq:num2}
\left|z\right| < \rho, \ 
\left|1-z\right| < \rho, \ 
\left|\frac{1}{z}\right| < \rho, \ 
\left|\frac{z-1}{z}\right| < \rho,
\left|\frac{1}{1-z}\right| < \rho, \ 
\left|\frac{z}{z-1}\right| < \rho, \ 
\end{equation}
with $0<\rho<1$.

\begin{figure}
\begin{center}
\epsfxsize=8cm \epsfbox{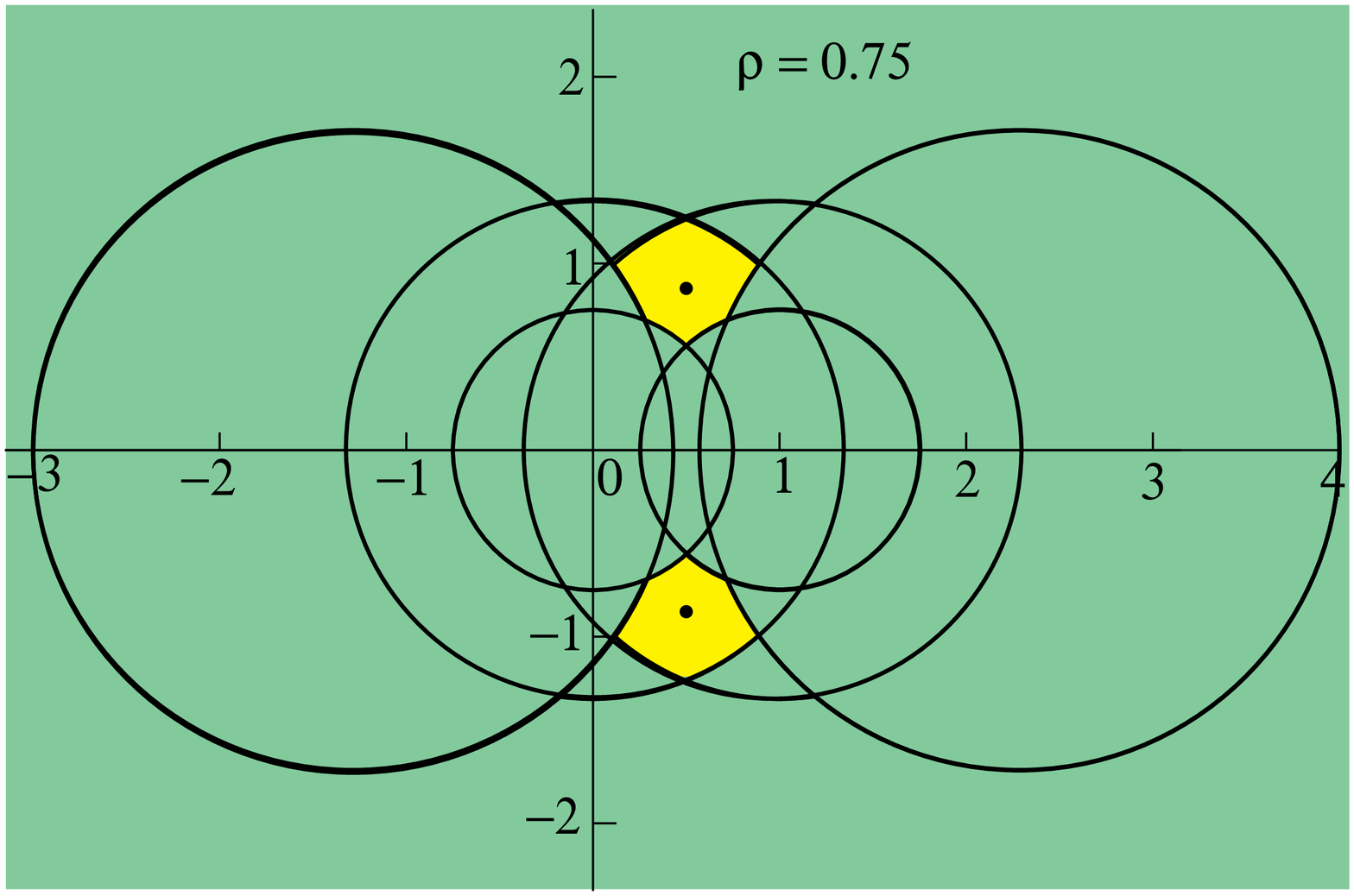}
\end{center}
\caption{\small In the light domains around the points
$e^{\pm \pi i/3}$ none of the inequalities of \pref{eq:num2} is satisfied.
\label{fig:abcn}}
\end{figure}

Not all points in the $z-$plane satisfy one of these inequalities for a given 
number $\rho$. In Figure \ref{fig:abcn} we take $\rho=\frac34$. 
In the dark area at least one
of the above inequalities is satisfied. In the light areas,  ``around'' the points
$e^{\pm \pi i/3}$, none of these inequalities is satisfied.
By choosing $\rho$ closer to unity these light domains become smaller.

For certain combinations of the parameters $a$, $b$ and $c$, 
the connection formulas become numerically unstable. For example,
if $c=a+b$, the relation in \pref{eq:b2} is well-defined, although 
two gamma functions are infinite. By using a limiting procedure the 
value of $_2F_1(a,b;a+b;z)$  can be found. For $c$ close to $a+b$ numerical 
instabilities occur when using \pref{eq:b2}. See \cite{Forrey:1997:CHF}
for many examples and details.

Other instabilities in the evaluation 
of the power series \pref{eq:i1} may arise for large values of $a$ and $b$.

In \cite[p. 71]{Wimp:1984:CRR} an example is given how to compute
a Gauss function with argument $z=e^{\pi i/3}$, the point that
is excluded
 from the convergence domains shown in Figure \ref{fig:abcn}. Wimp
considers the computation of
\begin{equation}\label{eq:num3}
\FG{\tfrac23}1{\tfrac43}{e^{\pi i/3}}=
\frac{2\pi e^{\pi i/6}\Gamma(\frac13)}{9[\Gamma(\frac23)]^2},
\end{equation}
by using a Miller algorithm for the hypergeometric functions
\begin{equation}\label{eq:num4}
f_n=\FG{n+a}{n+b}{2n+c}{z},\quad n=0,1,2,\ldots\ .
\end{equation}
This recursion type was not initially included in the group of 26 discussed in the
present paper. However as, discussed in \S \ref{sec:k2}, this recurrence can be related to the case
$k=2$ to conclude that $f_n$ is minimal. Therefore, Miller's algorithm can be applied
when a sum rule is provided, as done in \cite[p. 71]{Wimp:1984:CRR}.

We can also use basic form $k=13$, with pure $c-$recursion. From
\S\ref{sec:dk13} and \S\ref{sec:ak13} it follows that for the 
point $z=e^{\pi i/3}=\frac12+\frac12i\sqrt{3}$ the solutions
of the recursion relation \pref{eq:k131} are neither dominant nor minimal.
We use backward recursion for $_2F_1(a,b;c+n;z)$ with two starting values
for $n=29$ and $n=30$. With these large values of $c+n$ the power series
converges fast. In 15D arithmetic 
we have computed the value of \pref{eq:num3}
with a relative error $2\times 10^{-14}$. The exact value is
$$  0.883319375142724... + 0.509984679019064... i \ .$$
With recursion we obtain
$$   0.883319375142719 + 0.509984679019039 i \ .$$

In \cite[p. 72]{Wimp:1984:CRR} another example of the Miller algorithm is discussed for 
the basic form $k=13$ (pure $c-$recursion).

Several algorithms based on recursion relations for special cases of 
the Gauss functions have been published, in particular 
for computing Legendre functions. For recent papers, see 
\cite{Gil:1997:ELF},  \cite{Gil:2000:CTF} and \cite{Segura:1999:EAL}.

\section*{Acknowledgments}
A. Gil acknowledges financial support from Ministerio de Ciencia
y Tecnolog\'{\i}a (programa Ram\'on y Cajal). J. Segura acknowledges 
financial support from project BFM2003-06335-C03-02.
N.M. Temme acknowledges 
financial support from Ministerio de Educaci\'on y Ciencia
(Programa de Sab\'aticos) from project SAB2003-0113.

\bibliographystyle{plain}


\end{document}